\documentclass[a4paper,11pt]{article}

\usepackage{a4wide}
\usepackage{amsmath}
\usepackage{amssymb}
\usepackage{amsthm}
\usepackage{enumerate}
\usepackage[all]{xy}

\newtheoremstyle{custom}{9pt}{9pt}{\upshape}{}{\bfseries}{}{0.4em }{}

\swapnumbers \theoremstyle{plain}

\newtheorem{prop}[subsection]{Proposition}
\newtheorem{thm}[subsection]{Theorem}
\newtheorem{cor}[subsection]{Corollary}
\theoremstyle{definition}

\newtheorem{ex}[subsection]{Example}
\newtheorem{rem}[subsection]{Remark}
\theoremstyle{custom}
\newtheorem{nr}[subsection]{}

\newcommand{\df}{\emph}
\newcommand{\ev}{\mathrm{ev}}

\newcommand{\Id}{\mathrm{Id}}

\renewcommand{\c}{\mathsf{c}}
\newcommand{\cont}{\mathbf{cont}}
\newcommand{\op}{\mathrm{op}}
\newcommand{\tto}{\nrightarrow}
\newcommand{\two}{\mathbf{2}}

\newcommand{\f}{\mathfrak{f}}
\newcommand{\g}{\mathfrak{g}}

\newcommand{\x}{\mathfrak{x}}
\newcommand{\y}{\mathfrak{y}}
\newcommand{\z}{\mathfrak{z}}

\newcommand{\C}{\mathbf{C}}
\newcommand{\Alg}{\mathbf{Alg}}
\newcommand{\App}{\mathbf{App}}
\newcommand{\Cls}{\mathbf{Cls}}
\newcommand{\KlAlg}{\mathbf{KlAlg}}
\newcommand{\Kl}{\mathbf{Kl}}
\newcommand{\Mat}{\mathbf{Mat}}
\newcommand{\Ord}{\mathbf{Ord}}
\newcommand{\Set}{\mathbf{Set}}
\newcommand{\Sup}{\mathbf{Sup}}
\newcommand{\Top}{\mathbf{Top}}

\newcommand{\DM}{D_{\!_M}}
\newcommand{\FL}{F_{\!_\L}}
\newcommand{\FM}{F_{\!_M}}

\newcommand{\PM}{P_{\!_M}}

\newcommand{\TM}{T_{\!_M}}

\newcommand{\A}{\mathcal{A}}
\newcommand{\B}{\mathcal{B}}
\newcommand{\D}{\mathsf{D}}
\newcommand{\FF}{\mathfrak{F}}
\newcommand{\F}{\mathsf{F}}

\renewcommand{\L}{\mathbf{L}}
\renewcommand{\P}{\mathsf{P}}
\newcommand{\R}{\overline{\mathbf{R}}_+}

\newcommand{\T}{\mathsf{T}}

\newcommand{\V}{\mathbf{V}}
\newcommand{\X}{\mathfrak{X}}

\begin{document}

\title{A Kleisli-based approach to lax algebras}

\author{Gavin J. Seal\thanks{Financial support by the Swiss National Science Foundation is gratefully acknowledged.}}

\date{August 2006}

%Mathematics Subject Classification: 18C20, 18B30, 54A05

%Keywords: lax algebra, Kleisli category, topological category, topological space, fuzzy topological space, approach space, closure space

\maketitle

\begin{abstract}
By exploiting the description of topological spaces by either neighborhood systems or filter convergence, we obtain a neighborhood-like presentation of categories of lax algebras. A notable advantage of this approach is that it does not require the introduction of a lax extension of the associated monad functor. As a byproduct, the different philosophies underlying the construction of fuzzy topological spaces on one hand, and approach spaces on the other, may be simply expressed in terms of lax algebras.
\end{abstract}

%%%%%%%%%%%%%%%%%%%%%%%%%%%%%%%%%%%%%%%%%%%%%%%%%%%%%%%%%%%%%%%%%%%%%%%%%%%%%%%%%%%%%%%%%%%%%%%%%%%%%%%%%%%%%%%

\section{Introduction}

In \cite{Gahler:1992}, G\"ahler gave a presentation of the category $\Top$ of topological spaces as a category, denoted here by $\KlAlg(\F)$, of structured objects in the Kleisli category of the filter monad~$\F$:
\[
\Top\cong\KlAlg(\F)\ .
\]
This result lead to a natural definition of fuzzy topological spaces by extending the previous monad to a fuzzy filter monad. Lax algebras on the other hand (see \cite{Barr:1970}, \cite{Clementino/all:2004}, \cite{Hofmann/Tholen:200?} and \cite{Seal:2005}) provide a setting for the presentation of topological spaces as structured objects in the category of sets and relations:
\[
\Top\cong\Alg(\F,\two)\ .
\]
The category $\Alg(\F,\two)$ of these structured objects depends on the filter monad $\F$ and the two-element ordered chain $\two$. In this context, a notion of fuzziness may be introduced by replacing the ordered chain $\two$ by a larger unital quantale $\V$.

Although the previous descriptions of topological spaces are both based on the filter monad, the first approach does not require the existence of a lax extension of the monad functor, an extension which is crucial for the second. This remark is at the origin of the present work, as it suggests that the information pertaining to the construction of a lax extension of an arbitrary functor $T$ may be extracted from the objects of the category $\KlAlg(\T)$ associated with a monad $\T=(T,e,m)$. Not only is this the case, but the resulting category $\Alg(\T,\two)$ of lax algebras is isomorphic to $\KlAlg(\T)$, thus generalizing the previous correspondence obtained for the filter monad. By modifying Zhang's tower extension construction \cite{Zhang:2000} to include unital quantales, we can moreover define a category $\KlAlg(\T,\V)$ such that
\[
\KlAlg(\T,\V)\cong\Alg(\T,\V)\ ,
\]
and for which $\KlAlg(\T)$ is the $\V=\two$ instance: $\KlAlg(\T,\two)\cong\KlAlg(\T)$.

The present paper is organized as follows. After establishing certain definitions in Section~\ref{intro}, we present a preliminary result that puts forth conditions a lax extension should satisfy in order to generalize the isomorphism between $\KlAlg(\F)$ and $\Alg(\F,\two)$. These conditions lead in Section~\ref{laxextension} to the actual construction of a lax extension $\TM$ of $T$ from quite a different perspective than in \cite{Barr:1970}, \cite{Clementino/Hofmann:2004}, \cite{Schubert:2006} or \cite{Seal:2005}. Such a lax extension yields a category $\Alg(\T,\V)$ of $(\T,\V)$-algebras, where $\V$ may be any unital quantale rather than just $\two$. In turn, this leads in Section~\ref{main} to the definition of the category $\KlAlg(\T,\V)$ of Kleisli $(\T,\V)$-algebras, which is a generalization of the category $\KlAlg(\T)$. Our main result then conveniently states that the two categories $\KlAlg(\T,\V)$ and $\Alg(\T,\V)$ are isomorphic, thus allowing for an ``extension-free'' description of certain categories of lax algebras. In fact, if we refer to the original example, the Kleisli construction defines lax algebras via their ``neighborhood systems''. This is illustrated in Section \ref{clos}, in which the category $\Cls$ of closure spaces is presented by way of such structures (this also provides a new description of $\Cls$ as a category of lax algebras). Finally, Section \ref{fuz} is a brief incursion into the realm of fuzzy topology, in which fuzzy topological spaces (as defined for example in \cite{Hohle:2001}, or \cite{Gahler:1992}) are shown to be particular instances of lax algebras. This example is simply the original isomorphism $\KlAlg(\F)\cong\Alg(\F,\two)$ in which the filter monad is replaced by a fuzzy filter monad.

It seems relevant now to mention another example arising in the context of lax algebras. Indeed, recall that the category $\App$ of approach spaces is isomorphic to $\Alg(\F,\R)$, where $\R$ is the extended real line. Although the categories of fuzzy topological and approach spaces are both generalizations of $\Top\cong\Alg(\F,\two)$, the first is obtained by extending the filter monad $\F$, while the second by extending the underlying quantale $\two$. These examples clearly illustrate the difference between the two perspectives mentioned in the opening paragraph.

%%%%%%%%%%%%%%%%%%%%%%%%%%%%%%%%%%%%%%%%%%%%%%%%%%%%%%%%%%%%%%%%%%%%%%%%%%%%%%%%%%%%%%%%%%%%%%%%%%%%%%%%%%%%%%%

\section{Motivating result} \label{intro}

Before stating our preliminary result, we present a number of definitions, and recall some useful properties of the structures we will be using. For more details on lax algebras, we refer to the articles mentioned in the Introduction.

\begin{nr}\label{coherent}
\textbf{Monads factoring through a category.} Let $\C$ be a subcategory of the category $\Ord$ of preordered sets. A $\Set$-monad $\T=(T,e,m)$ \df{factors through} $\C$ if there is a functor $S:\Set\to\C$ that composes with the forgetful functor to yield $T$, and such that $m_X:T^2X\to TX$ is the image of a morphism $m_X:STX\to SX$ of $\C$. To simplify notations, we will not distinguish between $SX$ and $TX$; for example, if $\T$ factors through the category $\Sup$ of complete lattices and sup-preserving maps, $\Set$-maps $Tf:TX\to TY$, as well as $m_X:T^2X\to TX$, will be considered as a sup-preserving maps between complete lattices.

The monad $\T$ factors \df{coherently} through $\C$ if for any $f,g\in\Set(X,TY)$, we have
\begin{equation}\tag{$\ast$}
f\le g\implies m_Y\cdot Tf\le m_Y\cdot Tg\ ,
\end{equation}
where $\Set(X,TY)$ is equipped with the preorder induced by $TY$:
\begin{equation}\tag{$\ast\ast$}
f\le g\iff\text{ for all }x\in X,\text{ we have }f(x)\le g(x)\ .
\end{equation}
The notion of a monad factoring coherently though $\Sup$ is similar in spirit to the ordered monads of \cite{Gahler:1992}.
\end{nr}

\begin{nr}
\textbf{Kleisli $\T$-algebras.} Let $\T=(T,e,m)$ be a $\Set$-monad factoring coherently through $\Ord$, and denote by $\Kl(\T)$ the associated Kleisli category. Recall that the \df{Kleisli composition} $\beta\circ\alpha:X\to TZ$ of $\alpha:X\to TY$ and $\beta:Y\to TZ$ is given by $m_Z\cdot T\beta\cdot\alpha$, and the identity morphism in $\Kl(\T)$ is $e_X:X\to TX$. By composing with $e_Y:Y\to TY$, a $\Set$-map $f:X\to Y$ becomes an element of $\Set(X,TY)=\Kl(\T)(X,Y)$, and we write $f^\sharp:=e_Y\cdot f$. Remark that the condition ($\ast$) of \ref{coherent} is equivalent to preservation of the preorder on $\Set(X,TY)$ in the first variable of the Kleisli composition, while preservation of this preorder in the second variable simply follows from ($\ast\ast$).

The category $\KlAlg(\T)$ of \df{Kleisli $\T$-algebras}, has as objects pairs $(X,\alpha)$ with $X$ a set and
$\alpha:X\to TX$ a \df{structure map} that is \df{extensive} and \df{idempotent}:
\begin{enumerate}[$(K_1)$]
\item $e_X\le\alpha$ ,
\item $\alpha\circ\alpha\le\alpha$ .
\end{enumerate}
Of course, in presence of the extensivity condition, idempotency may be expressed as an equality. Morphisms $f:(X,\alpha)\to(Y,\beta)$ are $\Set$-maps $f:X\to Y$ satisfying:
\begin{enumerate}[$(K_1)$]
\setcounter{enumi}{2}
\item $f^\sharp\circ\alpha\le\beta\circ f^\sharp$ ,
\end{enumerate}
and composing as in $\Set$.
\end{nr}

\begin{nr}
\textbf{Complete distributivity.} In this work, $\V$ will always denote a unital quantale with two-sided unit $k$, and we will assume that $\V$ is \df{non-trivial}, that is, $\bot\ne k$. It will often be useful to suppose that $\V$ is \df{completely distributive}, \textit{i.e.} that any $b\in\V$ may be obtained as
\[
b=\bigvee\{a\in\V\,|\,a\prec b\}\ ,
\]
where $a\prec b$ means that for any subset $S\subseteq\V$ with $b\le\bigvee S$, there exists $s\in S$
satisfying $a\le s$. The following properties follow from the definition of $\prec$:
\begin{enumerate}[i)]
\item $a\prec b$ implies $a\le b$ ;
\item $a\le a'\prec b'\le b$ implies $a\prec b$ ;
\item $a\prec\bigvee S$ implies there exists $s\in S$ with $a\prec s$ .
\end{enumerate}
\end{nr}

\begin{nr}
\textbf{Lax extensions.} Let $\V$ be a unital quantale, and denote by $\Mat(\V)$ the category of \df{$\V$-matrices} (or \df{$\V$-relations}). Recall that the objects of $\Mat(\V)$ are sets, morphisms $r:X\tto Y$ are maps $r:X\times Y\to\V$, and the transpose $r^\circ:Y\tto X$ of $r:X\tto Y$ is defined by $r^\circ(y,x)=r(x,y)$ for all $x\in X$, $y\in Y$. Composition of $r:X\tto Y$ and $s:Y\tto Z$ is given by
\[
s\cdot r(x,z)=\bigvee_{y\in Y}r(x,y)\otimes s(y,z)\ ,
\]
and the identity $1_X:X\tto X$ is defined by $1_X(x,y)=k$ if $x=y$ and $1_X(x,y)=\bot$ otherwise. There is also an order on the hom-sets of $\Mat(\V)$ induced by the order on $\V$. Finally, a $\Set$-map $f:X\to Y$ will be identified with the matrix $f:X\tto Y$ given by $f(x,y)=k$ if $f(x)=y$ and $f(x,y)=\bot$ otherwise. We point out that if $f:X\to Y$, $g:W\to Z$ are $\Set$-maps, and $s:Y\tto Z$ is a $\V$-matrix, then
\[
g^\circ\cdot s\cdot f(x,w)=s(f(x),g(w))\ ,
\]
for all $x\in X$, $w\in W$.

A \df{lax extension} of a $\Set$-functor $T$ is a map
\[
\TM:\Mat(\V)\to\Mat(\V)\quad,\qquad(r:X\tto Y)\mapsto(\TM r:TX\tto TY)
\]
which preserves the order on the hom-sets and satisfies
\begin{enumerate}[$(T_1)$]
\item $Tf\le\TM f$ and $(Tf)^\circ\le\TM f^\circ$ ,
\item $\TM s\cdot\TM r\le\TM(s\cdot r)$ ,
\end{enumerate}
for all $f:X\to Y$, $r:X\tto Y$ and $s:Y\tto Z$. An important consequence of these conditions is that if
$f:X\tto Y$ and $g:Y\tto Z$ come from $\Set$-maps, then
\[
\TM(s\cdot f)=\TM s\cdot\TM f=\TM s\cdot Tf\quad\text{ and }\quad\TM(g^\circ\cdot r)=\TM g^\circ\cdot\TM r=(Tg)^\circ\cdot\TM r\ .
\]
Finally, when $TX$ is an ordered set, we say that $\TM$ is \df{order-compatible} if for all $\x,\y\in
TX$, we have
\[
\x\le\y\iff k\le\TM 1_X(\x,\y)\ .
\]
\end{nr}

\begin{nr}
\textbf{Lax algebras.} Let $\T=(T,e,m)$ be a $\Set$-monad equipped with a lax extension $\TM$ of $T$. The category $\Alg(\T,\V)$ of \df{$(\T,\V)$-algebras}, also called \df{lax algebras}, has as objects pairs $(X,r)$, where $X$ is a set, and $r:TX\tto X$ a \df{structure $\V$-matrix} satisfying the \df{reflexivity} and \df{transitivity} laws:
\begin{enumerate}[$(L_1)$]
\item $1_X\le r\cdot e_X$ ,
\item $r\cdot\TM r\le r\cdot m_X$ .
\end{enumerate}
Morphisms $f:(X,r)\to(Y,s)$ are $\Set$-maps $f:X\to Y$ satisfying:
\begin{enumerate}[$(L_1)$]
\setcounter{enumi}{2}
\item $r\le f^\circ\cdot s\cdot T f$ ,
\end{enumerate}
and composing as in $\Set$. In the case where the lax extension $\TM$ is order-compatible, then it follows that the structure matrix of a lax algebra $(X,r)$ reverses the order on $TX$, \textit{i.e.} for $\x,\y\in TX$ and $z\in X$, we have
\[
\x\le\y\implies r(\y,z)\le r(\x,z)\ .
\]
\end{nr}

\begin{rem}
In \cite{Seal:2005}, it was noted that a lax extension $\TM$ of $T$ naturally defined an order on $TX$. The order described therein was the opposite of the order given above in the definition of order-compatibility, so the structure matrices of the associated lax algebras \emph{preserved} that order rather than reversing it. In both cases however, the order on $TX$ is chosen as the natural one (for instance if $\T$ is the filter monad $\F$, then $\f\le\g$ \emph{always} means that the filter $\f$ is finer than $\g$).
\end{rem}

\begin{nr}
\textbf{Continuous lax algebras.} Let $\T=(T,e,m)$ be a $\Set$-monad factoring through $\Sup$, provided with an order-compatible lax extension $\TM$. The $(\T,\V)$-algebra $(X,r)$ is said to be \df{continuous} if for all $y\in X$ and $\A\subseteq TX$, we have
\[
\bigwedge_{\x\in\A}r(\x,y)=r(\bigvee\A,y)\ .
\]
The full subcategory of $\Alg(\T,\V)$ whose objects are the continuous lax algebras is denoted by $\Alg_\cont(\T,\V)$.

For example, if $\TM$ is the op-canonical extension of either the filter or the powerset monad, then any lax algebra is continuous (this is a particular case of Proposition \ref{prop3}).
\end{nr}

\begin{nr}
\textbf{Kleisli $\T$-algebras and $(\T,\two)$-algebras.} The correspondence between Kleisli $\F$-al\-ge\-bras and $(\F,\two)$-algebras is given as a ``functional description of lax algebras'' in \cite{Hofmann/Tholen:200?}. In fact, the case where $\T$ is the filter monad provides an ideal setting in which the relation between Kleisli $\T$-algebras and $(\T,\two)$-algebras may be described. Indeed, recall that a topological space may be defined by two conditions on its neighborhood filters (the Kleisli presentation), or by two conditions on the ``convergence'' relation between filters and points (the lax algebra presentation). Given the first, one can obtain the second by stating that every filter finer than the neighborhood filter of a point converges to that point. Similarly, if the relation between filters and point is given, one can obtain the neighborhood filter of a point by taking the coarsest among all the filters that converge to that point. This correspondence is concretized in the following result.
\end{nr}

\begin{prop} \label{prop1}
Let $\T=(T,e,m)$ be a $\Set$-monad factoring coherently through $\Sup$ provided with an order-compatible lax extension $\TM$.
\begin{enumerate}[i)]
\item There is a concrete functor $F:\Alg_\cont(\T,\two)\to\KlAlg(\T)$ that associates to a structure matrix $r:TX\tto X$ the structure map $\alpha_r:X\to TX$ given by
\[
\alpha_r(y):=\bigvee\{\x\in TX\,|\,r(\x,y)=\top\}\ .
\]
\item Suppose that the extension $\TM$ satisfies
\[
\TM r(\X,\y)=\top\implies m_X(\X)\le m_X\cdot T\alpha_r(\y)\ ,
\]
for all continuous structure matrices $r:TX\tto X$, and elements $\X\in T^2X$, $\y\in TX$. Then the concrete functor $G:\KlAlg(\T)\to\Alg_\cont(\T,\two)$ that associates to a structure map $\alpha:X\to TX$ the structure matrix $r_\alpha:TX\tto X$ given by
\[
r_\alpha(\x,y)=\top\iff\x\le\alpha(y)\
\]
is inverse to $F$.
\end{enumerate}
\end{prop}

\begin{proof}
The proof of this statement is almost identical to the proof of Theorem \ref{thm2}.
\end{proof}

\begin{rem}
In order to obtain a better description of the lax extension $\TM$, one might be tempted to replace the previous condition ``$\,\TM r(\X,\y)=\top\implies m_X(\X)\le m_X\cdot T\alpha_r(\y)\,$'' by a more restrictive one such as ``$\,\TM r(\X,\y)=\top\iff\X\le T\alpha_r(\y)\,$''. Although the result would remain true, this last equivalence is unfortunately not satisfied by the usual lax extensions of the powerset and filter monads.
\end{rem}

%%%%%%%%%%%%%%%%%%%%%%%%%%%%%%%%%%%%%%%%%%%%%%%%%%%%%%%%%%%%%%%%%%%%%%%%%%%%%%%%%%%%%%%%%%%%%%%%%%%%%%%%%%%%%%%

\section{A lax extension of $T$} \label{laxextension}

\begin{nr}
\textbf{The Kleisli extension.} Let $\T=(T,e,m)$ be a $\Set$-monad factoring coherently through $\Sup$, and consider the powerset monad $\P=(P,d,n)$. The unique sup-preserving map $\eta_X:PX\to TX$ extending $e_X:X\to TX$ along $d_X:X\to PX$ is given by
\[
\eta_X(A)=\bigvee\{e_X(x)\,|\,x\in A\}\ ,
\]
and defines a natural transformation $\eta:P\to T$ satisfying $\eta\cdot d=e$. Moreover, since $T\eta\cdot eP=eT\cdot\eta$ by naturality of $e$, we have for any $\A\in P^2X$ that
\[
m_X\cdot T\eta_X\cdot\eta_{PX}(\A)=\bigvee_{A\in\A}m_X\cdot T\eta_X\cdot e_{PX}(A)=\bigvee_{A\in\A}\eta_X(A)=\eta_X\cdot n_X(\A)\ .
\]
Thus, $\eta:\P\to\T$ is in fact a monad morphism.

For a $\V$-matrix $r:X\tto Y$, let $\rho_r=(\rho_r^a:Y\to PX)_{a\in\V}$ be the family of maps given by
\[
\rho_r^a(y):=\{x\in X\,|\,a\le r(x,y)\}\ .
\]
The \df{Kleisli extension of $T$} is $\TM:\Mat(\V)\to\Mat(\V)$ defined by
\[
\TM r(\x,\y):=\bigvee\{a\in\V\,|\,\x\le m_X\cdot T(\eta_X\cdot\rho_r^a)(\y)\}\ ,
\]
for all $\x\in TX$, and $\y\in TY$.
\end{nr}

\begin{prop}
The Kleisli extension $\TM$ of $T$ is a lax extension of $T$.
\end{prop}

\begin{proof}
Remark first that $\TM$ preserves the order on the hom-sets because $\T$ factors coherently through $\Sup$. For a map $f:X\to Y$, $a\in\V$, $x\in X$, and $y\in Y$, we have
\[
\rho_f^a(y)=\left\{\begin{array}{l@{\quad}l}
X & \text{if }a=\bot\\
f^{-1}\{y\} & \text{if }\bot\ne a\le k\\
\emptyset & \text{otherwise,}
\end{array}\right.
\qquad\text{and}\qquad \rho_{f^\circ}^a(x)=\left\{\begin{array}{l@{\quad}l}
Y & \text{if }a=\bot\\
\{f(x)\} & \text{if }\bot\ne a\le k\\
\emptyset & \text{otherwise.}
\end{array}\right.
\]
On one hand, for all $a\le k$ we have $\eta_X\cdot d_X\le\eta_X\cdot\rho_f^a\cdot f$, so $\x=m_X\cdot T(\eta_X\cdot d_X)(\x)\le m_X\cdot T(\eta_X\cdot\rho_f^a)\cdot Tf(\x)$; therefore, $k\le\TM f(\x,Tf(\x))$, and $Tf\le\TM f$. On the other hand, for all $a\le k$ we have $\eta_X\cdot d_X\cdot f\le\eta_X\cdot\rho_{f^\circ}^a$, and we may proceed as before to get $k\le\TM f^\circ(Tf(\x),\x)$, or $(Tf)^\circ\le\TM f^\circ$.

To prove that $\TM$ is a lax functor, let $r:X\tto Y$ and $s:Y\tto Z$ be two $\V$-matrices, $\x\in TX$, $\y\in
TY$, and $\z\in TZ$. Let $a,b\in\V$ be such that $\x\le m_X\cdot T(\eta_X\cdot\rho_r^a)(\y)$ and $\y\le m_Y\cdot T(\eta_Y\cdot\rho_s^b)(\z)$, so
$\x\le m_X\cdot T(\eta_X\cdot\rho_r^a)\cdot m_Y\cdot T(\eta_Y\cdot\rho_s^b)(\z)$. Note furthermore that
\begin{align*}
m_X\cdot T(\eta_X\cdot\rho_r^a)\cdot m_Y\cdot T(\eta_Y\cdot\rho_s^b)&=m_X\cdot m_{TX}\cdot T(T(\eta_X\cdot\rho_r^a)\cdot\eta_Y\cdot\rho_s^b)\\
&=m_X\cdot T(m_X\cdot\eta_{TX}\cdot P(\eta_X\cdot\rho_r^a)\cdot\rho_s^b)\\
&=m_X\cdot T(\eta_X\cdot n_X\cdot P\rho_r^a\cdot\rho_s^b)\ .
\end{align*}
Since
\[
n_X\cdot P\rho_r^a\cdot\rho_s^b(z)=\{x\in X\,|\,\exists y\in Y:a\le r(x,y),b\le
s(y,z)\}\subseteq\rho_{s\cdot r}^{a\otimes b}(z)\ ,
\]
we obtain $\x\le m_X\cdot T(\eta_X\cdot\rho_{s\cdot r}^{a\otimes b})(\z)$. Finally, suprema are preserved by $\otimes$ in each variable, so that $\TM r(\x,\y)\otimes\TM s(\y,\z)\le\TM(s\cdot r)(\x,\z)$, and $\TM$ is a lax extension of $T$ as claimed.
\end{proof}

\begin{prop} \label{prop3}
If $\TM$ is the Kleisli extension of $T$, then the following assertions hold.
\begin{enumerate} [i)]
\item If $\V$ is non-trivial, then $\TM$ is order-compatible.
\item For any set $X$ and $a\in\V$, $(a_{TX}\wedge 1_{TX})\le\TM(a_X\wedge 1_X)$, where $a_X:X\tto X$ is
the $\V$-matrix with constant value $a\in\V$.
\item Any $(\T,\V)$-algebra structure $r:TX\tto X$ is continuous, so that
\[
\Alg_\cont(\T,\V)=\Alg(\T,\V)\ .
\]

\end{enumerate}
\end{prop}

\begin{proof}
\textit{i)} Consider the identity $1_X:X\to X$. Then $\rho_{1_X}^a$ is the unit $d_X$ of the powerset monad whenever $a\in\V$ satisfies $\bot\ne a\le k$. Thus, on one hand $\x\le\y$ implies $k\le\TM 1_X(\x,\y)$. On the other hand, since $\V$ is non-trivial, $k\le\TM 1_X(\x,\y)$ implies there exists $a\in\V$ such that $\bot\ne a\le k$ and $\x\le m_X\cdot(T\eta_X\cdot\rho_{1_X}^a)(\y)=\y$.

\textit{ii)} Let $a\in\V\setminus\{\bot\}$ and consider the relation $r=a_X\wedge 1_X:X\tto X$. Then $\rho_r^b=d_X$ if in particular $b=a\wedge k$, so that $m_X\cdot(T\eta_X\cdot d_X)(\x)=\x$ yields $a\wedge k\le\TM r(\x,\x)$ as required.

\textit{iii)} Suppose now that $r:TX\tto X$ is a $(\T,\V)$-algebra structure, and let $\A\subseteq TX$. Then $\X=\eta_{TX}(\A)$ naturally satisfies $m_X(\X)=\bigvee\A$, so $\TM r(\X,e_X(z))\le r(\bigvee\A,z)$. By using naturality of $e$ and the definition of $\eta$, we observe that
\[
\TM r(\X,e_X(z))=\bigvee\big\{a\in\V\,\big|\,\textstyle{\bigvee}\{e_{TX}(\x)\,|\,\x\in\A\} \le\textstyle{\bigvee}\{e_{TX}(\x)\,|\,\x\in\rho_r^a(z)\}\big\}\ .
\]
Thus, for any $a$ such that $\A\subseteq\rho_r^a(z)$, we have $a\le \TM r(\X,e_X(z))$. This is the case in particular for $a=\bigwedge_{\x\in\A}r(\x,z)$, so that $\bigwedge_{\x\in\A}r(\x,z)\le r(\bigvee\A,z)$. We can conclude that this last inequality is in fact an equality since $r$ reverses the order in the first variable. \end{proof}

\begin{ex}
If $\V$ is a completely distributive lattice, then the Kleisli extensions of the powerset and the filter functors are given by
\begin{align*}
\PM r(A,B)&=\bigwedge_{x\in A}\bigvee_{y\in B}r(x,y)\ ,\qquad\text{ and }\\
\FM r(\f,\g)&=\bigwedge_{B\in\g}\,\bigvee_{A\in\f}\,\bigwedge_{x\in A}\,\bigvee_{y\in B}r(x,y)\
\end{align*}
respectively, where $r:X\tto Y$ is a $\V$-matrix, $A,B\in PX$, and $\f,\g\in FX$. Note that these are also the op-canonical extensions of the corresponding functors.
\end{ex}

%%%%%%%%%%%%%%%%%%%%%%%%%%%%%%%%%%%%%%%%%%%%%%%%%%%%%%%%%%%%%%%%%%%%%%%%%%%%%%%%%%%%%%%%%%%%%%%%%%%%%%%%%%%%%%%

\section{Towers of Kleisli algebras} \label{main}

In Proposition \ref{prop1}, it has been shown how Kleisli $\T$-algebras may be related to $(\T,\two)$-algebras. In order to extend this correspondence to other quantales than $\two$, we introduce the following definition, which is based on Zhang's tower extensions \cite{Zhang:2000}.

\begin{nr}
\textbf{Kleisli $(\T,\V)$-algebras.} Let $\T=(T,e,m)$ be a $\Set$-monad factoring coherently through $\Sup$, $\V$ a unital quantale, and notice that $\two$ embeds into $\V$ via
\[
\bot\mapsto\bot\ ,\quad\top\mapsto k\ .
\]
The \df{tower extension} of $\Kl(\T)$ along $\two\to\V$ is the category  $\KlAlg(\T,\V)$ (also called the category of \df{Kleisli $(\T,\V)$-algebras}) whose objects are pairs $(X,\alpha)$, with $\alpha$ a $\V$-indexed family of morphisms $\alpha=(\alpha^a:X\to TX)_{a\in\V}$ satisfying the following conditions:
\begin{enumerate}
\item[$(K_0)$] $\alpha^{\bigvee\A}=\bigwedge_{a\in\A}\alpha^a$ ,
\item[$(K_1)$] $e_X\le\alpha^k$ ,
\item[$(K_2)$] $\alpha^a\circ\alpha^b\le\alpha^{a\otimes b}$ ,
\end{enumerate}
for all $\A\subseteq\V$, and $a,b\in\V$. When the monad $\T$ is clearly determined by the context, such a structure $\alpha$ will be called a \df{$\V$-tower} on $X$. Morphisms $f:(X,\alpha)\to(Y,\beta)$ are maps $f:X\to Y$ such that
\begin{enumerate}
\item[$(K_3)$] $f^\sharp\circ\alpha^a\le\beta^a\circ f^\sharp$ ,
\end{enumerate}
for all $a\in\V$, and composing as in $\Set$. Since $(K_0)$ yields in particular that $\alpha^\bot(x)=\top$ for all $x\in X$, and tower $\alpha$, the category $\KlAlg(\T,\two)$ is concretely isomorphic to $\KlAlg(\T)$. Moreover, if $\V$ is completely distributive, then $(K_0)$ is equivalent to
\begin{enumerate}
\item[$(K_0')$] $\alpha^a=\bigwedge_{b\prec a}\alpha^b$ ,
\end{enumerate}
for all $a\in\V$. Notice that a $\V$-tower $\alpha$ is in fact a sup-preserving map $\alpha:\V\to\Set(X,TX)^\op$ that forms an op-lax functor with respect to the multiplicative structures. The previous presentation via families of Kleisli endomorphisms appears however to be more practical for our purpose.
\end{nr}

\begin{rem}
The original definitions of tower extensions in \cite{Brock/Kent:1997} and \cite{Zhang:2000} only considered the indexing set $\V$ as a complete lattice, rather than a quantale. This explains in part why approach spaces---which explicitly make use of the addition of $\R$ in their definition---were not directly described as tower extensions of topological spaces.
\end{rem}

\begin{thm} \label{thm2}
Let $(\T,\V)$ be a $\Set$-monad factoring coherently through $\Sup$, and suppose that $\V$ is completely
distributive. If $\Alg(\T,\V)$ denotes the category of lax algebras associated to the Kleisli extension $\TM$,
then $\Alg(\T,\V)$ and $\KlAlg(\T,\V)$ are concretely isomorphic.

More precisely, to a $(\T,\V)$-algebra structure $r:TX\tto X$ can be associated a $\V$-tower $\alpha_r=(\alpha_r^a:X\to TX)_{a\in\V}$ defined by
\[
\alpha_r^a(y):=\bigvee\{\x\,|\,a\le r(\x,y)\}\ ,
\]
and to a $\V$-tower $\alpha=(\alpha^a:X\to TX)_{a\in\V}$ can be associated a $(\T,\V)$-algebra structure $r_\alpha:TX\tto X$ given by
\[
r_\alpha(\x,y):=\bigvee\{a\in\V\,|\,\x\le\alpha^a(y)\}\ .
\]
This correspondence yields two concrete functors
\[
F:\Alg(\T,\V)\to\KlAlg(\T,\V)\quad\text{ and }\quad G:\KlAlg(\T,\V)\to\Alg(\T,\V)
\]
that are inverses of each other.
\end{thm}

\begin{proof}
To prove that $F$ is well-defined, consider a $(\T,\V)$-algebra structure $r:TX\tto X$. In order to verify $(K_0)$ for $\alpha_r$, let $y\in X$, and recall that $\rho_r^b(y)=\{\x\in TX\,|\,b\le r(\x,y)\}$ (where $b\in\V$), so we have $\alpha_r^b(y)=\bigvee\rho_r^b(y)$. If $\A\subseteq\V$, then by continuity of $r$, any $c\in\A$ satisfies
\[
c\le r(\textstyle{\bigvee}\rho_r^c(y),y)\le r(\textstyle{\bigwedge}_{b\in\A}\textstyle{\bigvee}\rho_r^b(y),y)\ ,
\]
so that $a=\bigvee_{c\in\A}c\le r(\bigwedge_{b\in\A}\alpha_r^b(y),y)$, and $\bigwedge_{b\in\A}\alpha_r^b(y)\le\alpha_r^a(y)$; equality follows, since $a\le b$ yields $\alpha_r^b\le\alpha_r^a$ for all $a,b\in\V$.

Reflexivity of $r$ immediately implies $(K_1)$. To verify $(K_2)$ it suffices to show that $a\otimes b\le r(\alpha_r^a\circ\alpha_r^b(y),y)$ for all $a,b\in\V$. For this, remark first that $r\cdot\alpha_r^b\ge(b_X\wedge 1_X)$, so $\TM(r\cdot\alpha_r^a)\ge(a_{TX}\wedge 1_{TX})$ by Proposition \ref{prop3}. Therefore,
\begin{align*}
r(m_X\cdot T\alpha_r^a\cdot\alpha_r^b(y),y)
&\ge\TM r(T\alpha_r^a\cdot\alpha_r^b(y),\alpha_r^b(y))\otimes r(\alpha_r^b(y),y)\\
&=\TM(r\cdot\alpha_r^a)(\alpha_r^b(y),\alpha_r^b(y))\otimes(r\cdot\alpha_r^b)(y,y)\ge a\otimes b
\end{align*}
by transitivity of $r$. Finally, if $f:(X,r)\to(Y,s)$ is a morphism of lax algebras, then $r(\x,y)\le s(Tf(\x),f(y))$ implies that $Tf\cdot\alpha_r^a(y)\le\alpha_s^a\cdot f(y)$ for all $a\in\V$.

Consider now a $\V$-tower $\alpha=(\alpha^a:X\to TX)_{a\in\V}$. To verify that $G$ is well-defined, we first need to prove the equality $\alpha_{r_\alpha}=\alpha$. For this, let $a\in\V$, $y\in X$, and set $\A=\{\x\in TX\,|\,a\le\bigvee\B_\x\}$, where $\B_\x=\{b\in\V\,|\,\x\le\alpha^b(y)\}$, so that $\alpha_{r_\alpha}^a(y)=\bigvee\A$. On one hand, $\x=\alpha^a(y)$ is in $\A$ (since $a$ is in $\B_\x$), so $\alpha\le\alpha_{r_\alpha}$. On the other hand, if $c\in\V$ is such that $c\prec a$, and $\x\in\A$, then by complete distributivity of $\V$ there exists $b\in\B_\x$ with $c\le b$. Therefore, for any $\x\in\A$ we have $\x\le\alpha^c(y)$, so $\x\le\bigwedge_{c\prec a}\alpha^c(y)=\alpha^a(y)$ by $(K_0')$. This implies $\alpha_{r_\alpha}\le\alpha$, as required.

Reflexivity of $r_\alpha$ is an immediate consequence of $(K_1)$. For transitivity, note that if $b\in\V$, then
\[
m_X\cdot m_{TX}\cdot T(\eta_{TX}\cdot\rho_{r_\alpha}^a)=m_X\cdot T(m_X\cdot\eta_{TX}\cdot\rho_{r_\alpha}^a)=m_X\cdot T\alpha_{r_\alpha}^a\ ,
\]
so that $\TM r_\alpha(\X,\y)\le\bigvee\{a\in\V,|\,m_X(\X)\le m_X\cdot T\alpha^a(\y)\}$ because $\alpha_{r_\alpha}\le\alpha$. By definition of $r_\alpha(\y,z)$, we have
\[
\TM r_\alpha(\X,\y)\otimes r_\alpha(\y,z)\le\bigvee\{a\otimes b\,|\,m_X(\X)\le\alpha^a\circ\alpha^b(z)\}\le
r_\alpha(m_X(\X),z)
\]
by $(K_2)$. Now, let $f:(X,\alpha)\to(Y,\beta)$ be a morphism, and suppose that $Tf\cdot\alpha^a(y)\le\beta^a\cdot f(y)$ for all $a\in\V$; this implies that $\{a\in\V\,|\,\x\le\alpha^a(y)\}\subseteq\{a\in\V\,|\,Tf(\x)\le\beta^a\cdot f(y)\}$, and consequently
$r_\alpha(\x,y)\le r_\beta(Tf(\x),f(y))$.

The proof that $r_{\alpha_r}=r$ is quite similar to that of $\alpha_{r_\alpha}=\alpha$. Indeed, let $\x\in TX$, and $y\in X$. Consider $\A=\{a\in\V\,|\,\x\le\bigvee\B_a\}$ where $\B_a=\{\y\in TX\,|\,a\le r(\y,y)\}$, so that $r_{\alpha_r}(\x,y)=\bigvee\A$. On one hand, we observe that $a=r(\x,y)$ is in $\A$ (since $\x$ is in $\B_a$), so $r\le r_{\alpha_r}$. On the other hand, if $a\in\V$ is such that $\x\le\bigvee\B_a$, then by continuity of $r$ we have that $a\le r(\bigvee\B_a,y)\le r(\x,y)$, so $r_{\alpha_r}\le r$. This shows that $GF=\Id$. Since $\alpha_{r_\alpha}=\alpha$ implies that $FG=\Id$, we conclude that $G$ is an isomorphism with inverse $F$.
\end{proof}

\begin{rem}
The previous theorem yields another presentation of the category $\App\cong\Alg(\F,\R)$ of approach spaces  \cite{Lowen:1997}, and suggests a new notion of ``approach system of neighborhoods''.
\end{rem}

%%%%%%%%%%%%%%%%%%%%%%%%%%%%%%%%%%%%%%%%%%%%%%%%%%%%%%%%%%%%%%%%%%%%%%%%%%%%%%%%%%%%%%%%%%%%%%%%%%%%%%%%%%%%%%%

\section{$\V$-valued closure spaces} \label{clos}

In \cite{Seal:2005}, it was shown that the category $\Cls$ of closure spaces could be seen as a category of $(\P,\two)$-algebras. However, since $\Top\cong\Alg(\F,\two)$ is a full subcategory of $\Cls$, it would be useful to describe $\Cls$ as a lax algebra of the form $\Alg(\D,\two)$, with the filter monad $\F$ appearing as a submonad of $\D$. The aim of this Section is to provide such a description.

\begin{nr}
\textbf{The up-set monad.} An \df{up-set} $\x$ on $X$ is a set of subsets of $X$ such that for any $A,B\subseteq X$, we have
\[
A\subseteq B\text{ and }A\in\x\implies B\in\x\ .
\]
The set $DX$ of up-sets on $X$ is equipped with the order given by reverse inclusion:
\[
\x\le\y\iff\x\supseteq\y\ .
\]
In fact, $DX$ is a complete lattice, with supremum obtained via intersection, and infimum via union. There is only one up-set containing the empty set, namely the set $PX$ of subsets of $X$, which is also the bottom element of $DX$. On the other hand, the empty set is the top element of $DX$.

The \df{up-set functor} $D$ assigns to a set $X$ the set $DX$ of up-sets on $X$, and sends a map $f:X\to Y$ to $Df:DX\to DY$ defined by
\[
B\in Df(\x)\iff f^{-1}(B)\in\x\ ,
\]
where $\x\in DX$. The \df{up-set monad} $\D$ is the triple $(D,e,m)$, where $e:\Id\to D$ and $m:D^2\to
D$ are the natural transformations whose components at $X$ are given by
\[
A\in e_X(x)\iff x\in A\qquad\text{ and }\qquad A\in m_X(\X)\iff A^\sharp\in\X\ ,
\]
where $A^\sharp=\{\x\in DX\,|\,A\in\x\}$, $x\in X$ and $\X\in D^2X$. It follows immediately from this definition that the filter monad is a submonad of $\D$.

In view of Theorem \ref{thm2}, we point out that $\D$ factors coherently through $\Sup$.
\end{nr}

\begin{rem}
If $\V$ is completely distributive, then Theorem \ref{thm2} allows us to describe the category of $(\T,\V)$-algebras associated to the Kleisli extension $\TM$ of $T$, without having to actually compute $\TM$. However, it is not difficult to verify that in the present case the Kleisli extension of $D$ is given by
\[
\DM r(\x,\y)=\bigwedge_{B\in\y}\bigvee_{A\in\x}\bigwedge_{x\in A}\bigvee_{y\in B}r(x,y) ,
\]
for a $\V$-matrix $r:X\tto Y$, and $\x,\y\in DX$.
\end{rem}

\begin{nr}
\textbf{$\V$-valued closure operators.} The objects of the category $\Cls(\V)$ are pairs $(X,c)$, where $X$ is a set and $c:PX\times X\to\V$ is a \df{$\V$-valued closure operator} (called a \df{closeness operator} in \cite{Seal:2005}), \textit{i.e.} a map satisfying:
\begin{enumerate}[$(C_1)$]
\item $x\in A\implies k\le c(A,x)$ ,
\item $A\subseteq B\implies c(A,x)\le c(B,x)$ ,
\item $a\otimes c(c_a[A],x)\le c(A,x)$ ,
\end{enumerate}
where $x\in X$, $A,B\subseteq X$, $a\in\V$ and $c_a[A]:=\{x\in X\,|\,a\le c(A,x)\}$. The pair $(X,c)$ is then a \df{$\V$-closure space}. A morphism of $\V$-closure spaces $f:(X,c)\to(Y,d)$ is a $\Set$-map $f:X\to Y$ satisfying $c(A,x)\le d(f(A),f(x))$. Recall that if $\V=\two$, then a closure operator $\gamma:PX\to PX$ may be defined via
\[
x\in\gamma(A)\iff c(A,x)=\top\ ;
\]
in fact, this equivalence yields a concrete isomorphism between $\Cls(\two)$ and the category $\Cls$ of closure spaces.
\end{nr}

\begin{nr}
\textbf{$\V$-graded closure operators.} Suppose that $\V$ is completely distributive. A \df{$\V$-graded closure operator} on $X$ is a family $\gamma=(\gamma^a:PX\to PX)_{a\in\V}$ of operators satisfying:
\begin{enumerate}[$(\Gamma_1)$]
\item $A\subseteq\gamma^a(A)$ for all $a\le k$ ,
\item $A\subseteq B$ implies $\gamma^a(A)\subseteq \gamma^a(B)$ ,
\item $\gamma^b\cdot\gamma^a(A)\subseteq\gamma^{a\otimes b}(A)$ ,
\item $\gamma^a(A)=\bigcap_{b\prec a}\gamma^b(A)$ ,
\end{enumerate}
for all $A,B\subseteq X$, and $a,b\in\V$. A morphism $f:(X,(\gamma^a)_{a\in\V})\to(Y,(\delta^a)_{a\in\V})$ is a $\Set$-map $f:X\to Y$ satisfying $f(\gamma^a(A))\subseteq\delta^a(f(A))$ for all $a\in \V$. For convenience, the pair $(X,(\gamma^a)_{a\in\V})$ is also called a \df{$\V$-closure space}, and the corresponding category denoted by $\Cls(\V)$ (this abuse is justified by the following proposition).
\end{nr}

\begin{prop}
If $\V$ is completely distributive, then the category of $\V$-closure spaces given by $\V$-valued closure operators is concretely isomorphic to the category of $\V$-closure spaces given by $\V$-graded closure operators.
\end{prop}

\begin{proof}
Suppose first that $(X,c)$ satisfies $(C_1)$ to $(C_3)$, and set $\gamma^a(A):=c_a[A]$ for $A\in PX$. Then $(C_1)$ clearly implies $(\Gamma_1)$, $(C_2)$ implies $(\Gamma_2)$, and it is not hard to see that $(C_3)$ implies $(\Gamma_3)$. For $(\Gamma_4)$, observe on one hand that $b\prec a$ implies $\gamma^a(A)\subseteq \gamma^b(A)$, so $\gamma^a(A)\subseteq\bigcap_{b\prec a}\gamma^b(A)$. On the other hand, if $x\in\bigcap_{b\prec a}\gamma^b(A)$, then $b\le c(A,x)$ for all $b\prec a$, so that $a\le c(A,x)$, as required. If $f:(X,c)\to(Y,d)$ satisfies $c(A,x)\le d(f(A),f(x))$ for all $A\subseteq X$, and $x\in X$, then $x\in\gamma^a(A)$ implies $f(x)\in\delta^a(f(A))$, where $\delta^a(B):=d_a[B]$.

For a pair $(X,(\gamma^a)_{a\in\V})$ satisfying $(\Gamma_1)$ to $(\Gamma_4)$, set $c(A,x):=\bigvee\{a\in\V\,|\,x\in\gamma^a(A)\}$. Then $(\Gamma_1)$ immediately implies $(C_1)$, and $(\Gamma_2)$ implies $(C_2)$. Let $a\in\V$, and remark that $c_a[A]=\gamma^a(A)$ by using complete distributivity of $\V$ and $(\Gamma_4)$. Thus, $a\otimes c(c_a[A],x)=\bigvee\{a\otimes b\,|\,x\in \gamma^b\cdot\gamma^a(A)\}\le c(A,x)$ by $(\Gamma_3)$. Finally, if $f:(X,(\gamma^a)_{a\in\V})\to(Y,(\delta^a)_{a\in\V})$ satisfies $f(\gamma^a(A))\subseteq \delta^a(f(A))$ for all $a\in \V$, then it follows that $c(A,x)\le d(f(A),f(x))$, where $d$ is the $\V$-valued closure operator assigned to $(\delta^a)_{a\in\V}$.

The fact that these correspondences are inverses of one another has been proved partially in the previous paragraph, and the remaining part is clear.
\end{proof}

\begin{prop} \label{prop56}
If $\V$ is completely distributive, then the category of Kleisli $(\D,\V)$-algebras is concretely isomorphic to the category $\Cls(\V)$ of $\V$-closure spaces.

More precisely, a $\V$-tower $\alpha=(\alpha^a: X\to DX)_{a\in\V}$ and a $\V$-graded closure operator $\gamma=(\gamma^a:PX\to PX)_{a\in\V}$ determine each other via
\[
x\in\gamma^a(A)\iff A^\c\notin\alpha^a(x)\ ,
\]
where $A\subseteq X$, and $A^\c$ denotes the complement of $A$ in $X$.
\end{prop}

\begin{proof}
Let $\alpha=(\alpha^a:X\to DX)_{a\in\V}$ be a $\V$-tower, and define $\gamma_\alpha=(\gamma_\alpha^a:PX\to PX)_{a\in\V}$ by $\gamma_\alpha^a(A):=\{x\in X\,|\,A^\c\notin\alpha^a(x)\}$, where $a\in\V$, and $A\subseteq X$. Since $\alpha$ satisfies $(K_0)$ and $(K_1)$, we have that $x\in A$ whenever $A\in\alpha^a(x)$ and $a\le k$. Thus, if $x\in A$, then $x\notin A^\c$, so $A^\c\notin\alpha^a(x)$ if $a\le k$, which proves $(\Gamma_1)$. If $A\subseteq B$ and $a\in\V$, $x\in X$ are such that $A^\c\notin\alpha^a(x)$, then $B^\c\subseteq A^\c$ implies that $B^\c\notin\alpha^a(x)$ because $\alpha^a(x)$ is an up-set, so we have $(\Gamma_2)$. To prove $(\Gamma_3)$, notice that for $a,b\in\V$ and $x\in X$,
\begin{align*}
A\in\alpha^a\circ\alpha^b(x)\iff A^\sharp\in D\alpha^a\cdot\alpha^b(x)\iff(\alpha^a)^{-1}(A^\sharp)\in\alpha^b(x)\ .
\end{align*}
Suppose now that $(\gamma_\alpha^a(A))^\c\notin\alpha^b(x)$. Since $(\gamma_\alpha^a(A))^\c=\{y\in X\,|\,A^\c\in\alpha^a(y)\}=(\alpha^a)^{-1}((A^\c)^\sharp)$, we have that $A^\c\notin\alpha^{a\otimes b}(x)$, as required. For $(\Gamma_4)$, we first note that if $b\prec a$, then $\alpha^a\le\alpha^b$, which implies that $\gamma_\alpha^a(A)\subseteq \gamma_\alpha^b(A)$ for all $A\subseteq X$. Thus, on one hand, we have $\gamma_\alpha^a(A)\subseteq\bigcap_{b\prec a}\gamma_\alpha^b(A)$. On the other hand, if $x\in\bigcap_{b\prec a}\gamma_\alpha^b(A)$, then $A^\c\notin\alpha^b(x)$ for all $b\prec a$. By $(K_0')$, we have $A^\c\notin\bigwedge_{b\prec a}\alpha^b(x)=\alpha^a(x)$, which proves $(\Gamma_4)$. If $f:(X,\alpha)\to(Y,\beta)$ satisfies $Df\cdot\alpha^a\le\beta^a\cdot f$ for all $a\in\V$, then $x\in \gamma_\alpha^a(A)\subseteq \gamma_\alpha^a(f^{-1}(f(A)))$, implies that $(f^{-1}(f(A)))^\c=f^{-1}(f(A)^\c)\notin\alpha^a(x)$, so $f(A)^\c\notin\beta^a(f(x))$, which shows that $f$ is a morphism of $\V$-closure spaces.

Suppose now that $\gamma=(\gamma^a:PX\to PX)_{a\in\V}$ is a $\V$-graded closure operator, and set $\alpha_\gamma^a(x):=\{A\in PX\,|\,x\notin \gamma^a(A^\c)\}$. To prove $(K_1)$, let $a\in\V$ be such that $a\le k$, and suppose that $A\in\alpha_\gamma^a(x)$. This implies that $x\notin \gamma^a(A^\c)$, so $x\notin A^\c$,
and $e_X(x)\le \alpha_\gamma^a(x)$. To verify $(K_2)$, recall from the previous paragraph that for $a,b\in\V$, we have $A\in\alpha_\gamma^a\circ\alpha_\gamma^b(x)$ if and only if $\{y\in X\,|\,A\in \alpha_\gamma^a(y)\}\in \alpha_\gamma^b(x)$. This last condition is equivalent to $\gamma^a(A^\c)^\c\in\alpha_\gamma^b(x)$, or $x\notin \gamma^b\cdot\gamma^a(A^\c)$, which yields $A\in\alpha_\gamma^{a\otimes b}(x)$, as required. Let now $f:(X,(\gamma^a)_{a\in\V})\to(Y,(\delta^a)_{a\in\V})$ be a morphism of $\V$-closure spaces. Then by using that $f(\gamma^a(f^{-1}(A^\c)))\subseteq\delta^a(f(f^{-1}(A^\c)))\subseteq\delta^a(A^\c)$, we observe that $A\in\alpha_\delta^a\cdot f(x)$ implies $f(x)\notin\delta^a(A^\c)$, so $x\notin \gamma^a(f^{-1}(A^\c))=\gamma^a(f^{-1}(A)^\c)$, and $f^{-1}(A)\in(\alpha_\gamma)^a(x)$. Therefore, $A\in\alpha_\delta^a(f(x))$ yields $A\in Df(\alpha_\gamma)^a(x)$, and $f$ is a morphism of Kleisli $(\T,\V)$-algebras.

To show that the previous correspondences yield an isomorphism, note that
\[
\gamma_{\alpha_\gamma}^a(A)=\{x\in X\,|\,x\in \gamma^a(A)\}=\gamma^a(A)\ ,
\]
for all $a\in\V$, and $A\subseteq X$. Furthermore,
\[
\alpha_{\gamma_\alpha}^a(x)=\{A\in PX\,|\,A\in\alpha^a(x)\}=\alpha^a(x)\ ,
\]
for all $a\in\V$, and $x\in X$, so we are done.
\end{proof}

\begin{rem}
In the case $\V=\two$, the Kleisli algebras present closure spaces by way of their ``neighborhood systems", where the neighborhood of a point $x\in X$ with respect to a closure operator $\gamma:PX\to PX$ is given by $\mathcal{N}(x)=\{A\subseteq X\,|\,x\in \gamma(A^\c)^\c\}$. Of course, this immediately leads us to the definition of an interior operator on a set $X$, and it is well known that the category of closure spaces may also be described by such operators.
\end{rem}

\begin{cor}
If $\V$ is completely distributive, then
\[
\Alg(\D,\V)\cong\Cls(\V)\ ,
\]
where $\Alg(\D,\V)$  is the category of $(\D,\V)$-algebras associated to the Kleisli extension of $D$.
\end{cor}

\begin{proof}
This is a direct consequence of the previous proposition and Theorem \ref{thm2}. Note that the complete distributivity condition on $\V$ allows us to treat a number of isomorphisms in one stroke, but it is not necessarily the most efficient hypothesis for each one. For example, Schubert proved a similar result (\cite{Schubert:2006}, 4.4.2) by only assuming that $k=\top$.
\end{proof}

\begin{rem}
In the case $\V=\two$, the corollary suggests that up-sets might play the role of filters in a convergence theory for closure spaces. However, other candidates appear in the literature, such as the $p$-stacks of \cite{Kent/Min:2002}, or the rasters of \cite{Giuli/Slapal:2005}. Nonetheless, in the present context all these concepts are very similar. In the case $\V=\two$ for example, whenever $\alpha:X\to DX$ satisfies $e_X\le\alpha$, then the intersection of all elements of $\alpha(x)$ is non-empty, so that $\alpha(x)$ is both a $p$-stack and a raster (as long as the empty set is also considered to be such a structure), and the structures of the Kleisli $\D$-algebras restrict accordingly.
\end{rem}

%%%%%%%%%%%%%%%%%%%%%%%%%%%%%%%%%%%%%%%%%%%%%%%%%%%%%%%%%%%%%%%%%%%%%%%%%%%%%%%%%%%%%%%%%%%%%%%%%%%%%%%%%%%%%%%

\section{Many-valued topologies} \label{fuz}

\begin{nr}
\textbf{The $\L$-valued filter monad.} (See \cite{Hohle:2001}) Let $\L$ be a complete lattice provided with a binary operation $\ast$ that is monotone in both variables (in particular, any quantale $\V$ is such a lattice with its binary operation given by $\otimes$; better yet, any complete lattice has a binary operation given by infimum). There is an induced order on the set $\L^X$ of maps from $X$ to $\L$ defined by
\[
A\le B\iff\text{ for all }x\in X,\text{ we have }A(x)\le B(x)\ ,
\]
where $A,B\in\L^X$. The top and bottom elements of $\L$ are denoted by $\top$ and $\bot$, and those of $\L^X$ by $\top^X$ and $\bot^X$ respectively. An $\L$-valued filter on $X$ is a map $\f:\L^X\to\L$ satisfying the following conditions for $A,B\in\L^X$:
\begin{enumerate}[$(F_1)$]
\item $\f(\top^X)=\top$ ,
\item $A\le B$ implies $\f(A)\le\f(B)$ ,
\item $\f(A)\ast\f(B)\le\f(A\ast B)$ .
\end{enumerate}
The set of all $\L$-valued filters on $X$ is denoted by $\FL X$. Of course, in the case $\L=\two$, we get the usual definition of filters, so that $F_{\!_\two}X=FX$.

The \df{$\L$-valued filter functor} $\FL$ assigns to a set $X$ the set $\FL X$, and sends a map $f:X\to Y$ to $\FL f:\FL X\to \FL Y$ defined by
\[
[\FL f(\f)](A)=\f(A\cdot f) ,
\]
where $A\in\L^Y$, and $\f\in\FL X$. The \df{$\L$-valued filter monad} $\F_{\!_\L}$ is the triple $(\FL,e,m)$, where $e:\Id\to\FL$ and $m:\FL^2\to\FL$ are the natural transformations whose components at $X$ are obtained via
\[
[e_X(x)](A)=A(x)\qquad\text{ and }\qquad [m_X(\FF)](A)=\FF(\ev_A)\ ,
\]
where $A\in\L^X$, $\FF\in\FL^2 X$, and $\ev_A:\FL X\to\L$ is given by $\ev_A(\f)=\f(A)$. The order on $\FL X$ is defined by
\[
\f\le\g\iff\g(A)\le\f(A)\text{ for all }A\in\L^X\ .
\]
With this order, $\FL$ is a complete lattice, and it is easily checked that $\FL$ factors coherently through $\Sup$ (see also \cite{Hohle:2001}, Proposition 2.4.2.3).

The category of Kleisli $F_{\!_\L}$-algebras is called the category of \df{$\L$-valued topological spaces} (or \df{fuzzy topological spaces}, see \cite{Gahler:1995b}), and is denoted by $\Top(\L)$. Furthermore, if $(X,\alpha)$ is a Kleisli $F_{\!_\L}$-algebra, then for each $x\in X$, the $\L$-valued filter $\alpha(x)$ is called the \df{$\L$-valued neighborhood system} of $x$.
\end{nr}

\begin{rem}
The previous definition of the $\L$-valued filter monad differs in two points with the one given in \cite{Hohle:2001}. First, we do not ask that filters $\f\in\FL X$ satisfy $\f(\bot^X)=\bot$. This will not be a problem in using results of \textit{op.cit.}, since $\L$-valued topologies are defined via neighborhoods, and a neighborhood $\f$ of $x$ must satisfy $\f(A)\le A(x)$ for all $A\in\L^X$; in particular $\f(\bot^X)=\bot$, so the resulting $\L$-valued neighborhood systems are the same. Second, the order on $\FL X$ is chosen as opposite to the one in the cited reference. In particular, with the present definition, $\FL X$ is a complete lattice rather than an ``almost complete join-semilattice''.
\end{rem}

\begin{prop} \label{propfuz}
The category of $(\F_{\!_\L},\V)$-algebras associated to the Kleisli extension of $\FL$ is isomorphic to the category of $\L$-valued topologies:
\[
\Alg(\F_{\!_\L},\two)\cong\Top(\L)\ .
\]
\end{prop}

\begin{proof}
This is a direct consequence of Propositions \ref{prop1} and \ref{prop3}.
\end{proof}

%%%%%%%%%%%%%%%%%%%%%%%%%%%%%%%%%%%%%%%%%%%%%%%%%%%%%%%%%%%%%%%%%%%%%%%%%%%%%%%%%%%%%%%%%%%%%%%%%%%%%%%%%%%%%%%

\textbf{Acknowledgements.} The author is particularly indebted to Christoph Schubert for pointing out in the first place that the Kleisli presentation of topological spaces did not require a lax extension of the monad functor, for mentioning Brock and Kent's limit tower spaces \cite{Brock/Kent:1997}, and in general for a number of insightful discussions and pertinent remarks. He also wishes to thank Walter Tholen for his suggestions towards improving the presentation of this paper.

\bibliographystyle{plain}

\end{document}